\documentclass[12pt,draftcls,onecolumn]{IEEEtran}

\newtheorem{TT}{Theorem}

\newtheorem{DD}{Definition}

\ifCLASSINFOpdf
\else
\fi

\ifCLASSOPTIONcompsoc
\usepackage[nocompress]{cite}
\else
\usepackage{cite}
\fi

\usepackage{amsmath}
\usepackage{amssymb}
\usepackage{graphicx}
\usepackage{epstopdf}
\usepackage{mathtools}
\usepackage{xcolor}

\begin{document}

\title{A Simple Introduction to Free Probability Theory and Its Application to Random Matrices}


\author{
        Xiang-Gen Xia, \IEEEmembership{Fellow}, \IEEEmembership{IEEE} 

\thanks{
X.-G. Xia is with the Department of Electrical and Computer Engineering, University of Delaware, Newark, DE 19716, USA (e-mail: xxia@ece.udel.edu).}

}

\maketitle


\begin{abstract}
Free probability theory started in the 1980s has attracted much attention lately in signal processing and communications areas due to its applications in 
large size 
random matrices. However, it involves with massive mathematical concepts
and notations, and is really hard for a general reader to comprehend. The main goal of this paper is to 
briefly describe this theory and its application in random matrices
 as simple as possible so that it is easy to follow. 
Applying 
free probability theory, one is able to 
calculate the distributions of the eigenvalues/singular-values of 
large size random matrices using only the second order statistics 
of the matrix entries. One of such applications is the mutual 
information calculation of a massive MIMO system.
\end{abstract}



\section{Introduction}\label{sec1}
Free probability theory was started by Voiculescu in the 1980's \cite{voi1, voi2, voi3}. It is about calculating moments (or distributions) of non-commutative random variables, such as, random matricies where the matrix entries are classical random variables. 

In classical probability theory, random variables are usually real-valued and can be extended to be complex-valued.
For convenience, let us say that they are real-valued. 
 Therefore, they are commutative.  For example, assume  $x_1, x_2$ are 
two independent non-zero random variables and $E$ denotes the expectation. Then,\begin{equation}\label{1.1} 
E(x_1x_2x_1x_2)=E(x_1^2x_2^2)=E(x_1^2)E(x_2^2)>0,
\end{equation}
 no matter whether $x_1$ and/or $x_2$ have $0$ mean or not, 
which is because $x_1$ and $x_2$ are commutative. 

However, if $x_1$ and $x_2$ are not commutative, then,
the property (\ref{1.1}) 
 may not hold and two natural questions are as follows.
 What will happen 
to (\ref{1.1})? What does the independence mean to non-commutative 
random variables? 

Free probability theory addresses the above two questions. 
It introduces freeness between non-commutative random variables,  which is 
analogous to the independence between classical commutative random variables. 
It basically says that although $E(x_1x_2x_1x_2)$ may not be equal to 
$E(x_1^2x_2^2)$, it is $0$ if $x_1$ and $x_2$ are free and  both 
have mean $0$. 

With this freeness, when a large number of free random variables are summed
with proper weights, it converges to the classical semicircular distribution. This is the free central limit theorem
similar to the classical central limit theorem,
where Gaussian distribution corresponds to semicircular distribution. 
Note that the eigenvalue distribution of a random matrix with entries of independent Gaussian random variables 
(for simplicity, the matrix symmetricity is not specified here) 
goes to semicircular distribution as well when the matrix size goes to infinity. This suggests a connection between free random variables and large size random matrices. 
Free probability theory 
says that, it indeed has a strong connection, i.e., random matrices of 
independent Gaussian random variables become free when the matrix size goes to infinity. In other words, when the size of matrices is large, these matrices 
are approximately free. 

Furthermore, the entries in random matrices can be 
replaced by free semicircular random variables (called deterministic equivalent). With the replacement,  all the joint moments or cumulants of random matrices 
can be calculated, which may lead to the calculations of the distributions of 
the eigenvalues of the 
functions of these random matrices. 

This is the reason why free probability theory has attracted much attention in wireless communications and signal processing areas. Massive MIMO systems have been identified as potential candidates in future wireless communications systems. In massive MIMO systems, their channel matrices are random of large sizes. Therefore, it is natural to apply free probability theory to do some of 
the difficult calculations, such as, channel capacity \cite{tul1, cou, anan}. 
It is particularly interesting when some statistics of a channel 
matrix of large size, such as, the first two moments (covariances)
 of the channel coefficients, are known, how we calculate the 
channel performance without 
performing Monte Carlo simulations that may be hard to do 
in practice when the channel matrix 
size is large,  such as, a massive MIMO channel.

The main goal of this tutorial paper is to briefly introduce free probability theory and its application to large size random matrices so that an ordinary   
researcher in signal processing and communications areas can easily understand. 

In the following, we adopt most of the notations in  
Speicher \cite{spe1,spe2,spe3,spe4}. 
All the results described below are from \cite{spe1,spe2,spe3,spe4} as well. 
The remainder of this paper is organized as follows. 
In Section \ref{sec2}, we describe the basics of free random variables and 
the free central limit theorem without proof. 
In Section \ref{sec3}, we describe the calculations/relations of joint moments, cumulants, and distributions of multiple free random variables.
In Section \ref{sec4}, we describe random matrices and the approximate
distributions of their eigenvalues. In Section \ref{sec5}, we
describe free deterministic equivalents for random matrices. We also describe
how to calculate  
the Cauchy transforms of random matrices using the  second order 
statistics of their entries. In Section \ref{sec6}, we conclude this 
paper.

\section{Free Random Variables}\label{sec2}
For convenience, in the following we will use as simple 
notations as possible, which may be too simplified in terms of mathematical 
rigorousness.

Let $x_1, x_2,..., x_n$ be $n$ elements that may not be commutative, and $E$ be a linear functional on these elements so that $E(1)=1$. 
Examples of these elements are matrices and $E$ is like the expectation of a classical random variable. 

\begin{DD}\label{def1}
Elements (or random variables) $x_1,x_2,..., x_n$ are called free or freely independent, if
for any $m$ polynomials $p_k(x)$, $1\leq k\leq m$, with $m\geq 2$, 
\begin{equation}\label{2.1}
E(p_1(x_{i_1})p_2(x_{i_2})\cdots p_m(x_{i_m}))=0,
\end{equation}
when $E(p_k(x_{i_k}))=0$ for all $k$, $1\leq k\leq m$, and any two neighboring 
indices $i_l$ and $i_{l+1}$ are not equal, i.e., 
$1\leq i_1\neq i_2 \neq  \cdots \neq i_m\leq n$.
\end{DD}

From (\ref{2.1}), if $x_1$ and $x_2$ are free, then 
$E(x_1x_2x_1x_2)=0$ when $E(x_1)=E(x_2)=0$, where $m=4$,
$i_1=1, i_2=2, i_3=1, i_4=2$,  and polynomials 
$p_k(x)=x$ for $1\leq k\leq 4$. Comparing with (\ref{1.1}) in the classical commutative case, 
independent real-valued random variables are not free. 
The terminology ``free'' comes from the concept of free groups, where there is 
no  any nontrivial relation between any generating elements of a free group. 

One might want to ask why, in the above definition, polynomials of the 
random variables $x_k$ are used. It is for the convenience later in calculating their joint moments. Note that 
in free probability theory context, it  not convenient to directly define
 density functions (or distribution functions) for noncommutative random variables. However, as we can recall, in the classical probability theory, if all the moments of a random variable are known, its characteristic function can be often determined and therefore, its density function can be often determined as well. Thus, calculating all the joint moments of free random variables may be sufficient for their joint distributions. Its details will be described in Section \ref{sec3}. 

The set ${\cal A}_k$ of all polynomials $p(x_k)$ of $x_k$ including the identity element $1=x_k^0$ is called the subalgebra generated by element $x_k$ for 
$1\leq k\leq n$. 
Subalgebras ${\cal A}_1, {\cal A}_2,...,{\cal A}_n$ are called free if and only if elements $x_1, x_2, ..., x_n$ are free. 
Clearly, when elements $x_1, x_2, ..., x_n$ are free, 
for any $n$ polynomials $p_1(x),...,p_n(x)$, elements 
$p_1(x_1),...,p_n(x_n)$ are free as well. 

If elements $x_1, x_2, \cdots, x_n$ are free, they are called free
 random variables. With the above freeness definition, although one may  
construct abstract free random variables using possibly many 
mathematical concepts, it is not easy to
show concrete examples of free random variables at this moment. 

Two sets ${\cal S}_1$ and ${\cal S}_2$ are called free if any element 
in ${\cal S}_1$ and any element in ${\cal S}_2$ are free. 
With property (\ref{2.1}), when $\{x_1,x_3\}$ and $x_2$ are free, 
it is easy to check 
that $E(x_1x_2)=E(x_1)E(x_2)$ and $E(x_1x_2x_3)=E(x_1x_3)E(x_2)$.

In many practical applications, we may need to deal with complex-valued random variables, such as, complex Gaussian, 
 where the complex conjugation $*$ is usually used. In correspondence 
with the complex conjugation, the above freeness becomes $*$-freeness.
We call that  $x_1,x_2,\cdots,x_n$ are $*$-free, if (\ref{2.1}) holds
when the polynomials $p_k(x)$ in Definition \ref{def1} are changed to 
polynomials $p_k(x, x*)$ of two variables. 
If $x=x^*$, element $x$ is called self-adjoint. For example, when $x$ is a 
matrix and $*$ is the complex conjugate transpose operation,  if $x$ is Hermitian, then $x$ is self-adjoint. In this case, $x$ can be diagonalized by a unitary matrix and all its eigenvalues are real-valued.

\begin{DD}\label{def2}
\begin{itemize}
\item[1)] When two random variables $x_1$ and $x_2$ have all the moments the same, 
i.e., $E(x_1^i)=E(x_2^i)$ for all positive integers $i$, 
 they are called identically distributed or  having the same distribution. 
\item[2)] For a sequence of random variables $x_n$, $n=1,2,...$, we call 
$x_n$ converges to $x$ in distribution when $n$ goes to infinity,
if all the moments of $x_n$ converge to the moments of $x$ as $n$ goes 
to infinity, i.e., for any positive integer $m$,
\begin{equation}\label{2.2}
\lim_{n\rightarrow \infty} E(x_n^m) =E(x^m),
\end{equation}
which is denoted as $\lim_{n\rightarrow \infty} x_n \overset{distr}{=} x$
or $x_n \overset{distr}{\longrightarrow} x$ as $n\rightarrow \infty$. 
\item[3)] Let $I$ be an index set. For each $i\in I$, let $x_n^{(i)}$, $n=1,2,...$, be a sequence of random variables. We call that $(x_n^{(i)})_{i\in I}$ converges to $(x^{(i)})_{i\in I}$ in distribution, 
if 
\begin{equation}\label{2.02}
\lim_{n\rightarrow \infty}E(x_n^{(i_1)}\cdots x_n^{(i_k)})=E(x^{(i_1)}\cdots x^{(i_k)})
\end{equation}
for all positive integers $k$ and all $i_1,..., i_k\in I$, 
which is denoted as 
$$
\lim_{n\rightarrow \infty} (x_n^{(i)})_{i\in I} \overset{distr}{=}
 (x^{(i)})_{i\in I} \mbox{ or }
(x_n^{(i)})_{i\in I} \overset{distr}{\longrightarrow}
 (x^{(i)})_{i\in I} \mbox{ as }n\rightarrow \infty.
$$
\end{itemize}
\end{DD}

The definition in 2) is about the convergence in distribution for 
a single sequence of random variables and the definition in 3)
is about the convergence in distribution for multiple sequences of random variables jointly. 

One of the most important results in classical probability theory is the central limit theorem. It says that the summation of independent random variables of
a totally fixed variance converges to Gaussian random variable, when the number of the independent random variables goes to infinity.  For free random variables, it has the following free central limit theorem.

\begin{TT}\label{thm1}
Let $x_k$, $k=1,2,...$, be a sequence of self-adjoint, freely independent, and identically distributed random variables  with $E(x_k)=0$ and $E(x_k^2)=\sigma^2$.  For a positive integer $n$, let 
\begin{equation}\label{2.3}
S_n=\frac{x_1+x_2+\cdots +x_n}{\sqrt{n}}.
\end{equation}
Then, $S_n$ converges in distribution to a semicircular element $s$ of variance $\sigma^2$ as $n\rightarrow \infty$, i.e.,
\begin{equation}\label{2.4}
\lim_{n\rightarrow \infty}E(S_n^{i})=
\left\{ \begin{array}{ll} \sigma^{i}C_{i/2}, & \mbox{if $i$ is even},\\
                          0, & \mbox{if $i$ is odd},
        \end{array}
\right.
\end{equation}
where $C_k$ is the Catalan number and the $(2k)$th moment of the semicircular 
distribution:
\begin{equation}\label{2.5}
C_k=\frac{1}{2\pi} \int_{-2}^{2} t^{2k} \sqrt{4-t^2} dt = \frac{1}{k+1}
{2k \choose k}. 
\end{equation}
\end{TT}

The random variable $s$ in Theorem \ref{thm1} is 
called a semicircular element in this context
and it, after divided by $\sigma$, has the same distribution as the classical semicircular 
random variable  of density function
\begin{equation}\label{2.6}
q(t)=\left\{ \begin{array}{ll} \frac{1}{2\pi} \sqrt{4-t^2},& \mbox{if }|t|<2,\\
                      0, & \mbox{otherwise}.
\end{array}\right.
\end{equation}
 Its mement of an even order has the form in (\ref{2.4}) and an odd order is always $0$. 

Note that semicircular distributions are the asymptotic distributions of the eigenvalues 
of Hermitian Gaussian random matrices when the matrix size goes to infinity, 
which is called Wigner's semi-circle law and will be discussed in more details in Section \ref{sec4} later.

\section{Moments, Cumulants, and Cauchy Transforms}\label{sec3}
As mentioned earlier, it is not convenient to directly 
define a density function or probability measure 
for a noncommutative random variable, and instead 
 its all moments are defined and the freeness is to simplify the joint moments 
between free random variables. 

In order to see how moments are related to 
distributions of free random variables, 
let us first see how in classical probability theory, a probability measure and its moments are related. 

Let $\mu(t)$ be a probability measure on the real line $\mathbb{R}$. Assume its all moments are finite and let $m_i$ be its $i$th moment for 
a positive integer $i$  and $\phi(t)$ be its characteristic function, i.e.,
\begin{equation}\label{3.1}
m_i=\int_{\mathbb R} t^i d\mu(t),   \mbox{ and  }
\phi(t)=\int_{\mathbb R} e^{{\bf i}\tau t}d\mu(\tau),
\end{equation}
where ${\bf i}\overset{\Delta}{=}\sqrt{-1}$. 
Then, it is easy to see 
\begin{equation}\label{3.2}
m_i={\bf i}^{-i}\phi^{(i)}(0), \mbox{ and }
\phi(t)=\sum_{i=0}^{\infty} m_i \frac{({\bf i}t)^i}{i!},
\end{equation}
where $\phi^{(i)}(t)$ stands for the $i$th derivative
of $\phi(t)$. Furthermore, we can write
\begin{equation}\label{3.3}
\log(\phi(t))=\sum_{i=1}^{\infty} k_i \frac{({\bf i}t)^i}{i!} \mbox{ with }
k_i = {\bf i}^{-i} \left. \frac{d^i}{dt^i} \log(\phi(t))\right|_{t=0},
\end{equation}
where $k_i$ are called the cumulants of $\mu(t)$. We will call them the 
classical cumulants. The moment sequence $\{m_i\}_{i\geq 0}$ and the cumulant sequence $\{k_i\}_{i\geq 1}$ can be determined 
from each other: 
\begin{eqnarray}
m_n & = & \mathop{\sum_{1\cdot r_1+\cdots + n\cdot r_n=n}}_{r_1,...,r_n\geq 0}
         \frac{n!}{(1!)^{r_1}\cdots (n!)^{r_n} r_1!\cdots r_n!} 
          k_1^{r_1} \cdots k_n^{r_n} \label{3.4}\\
k_n & = & \mathop{\sum_{1\cdot r_1+\cdots + n\cdot r_n=n}}_{r_1,...,r_n\geq 0}
         \frac{(-1)^{r_1+\cdots +r_n-1}(r_1+\cdots r_n-1)! n!}{(1!)^{r_1}\cdots (n!)^{r_n} r_1!\cdots r_n!} 
          m_1^{r_1} \cdots m_n^{r_n}. \label{3.5}
\end{eqnarray}
Sometimes, cumulants may be easier to obtain than moments. In this case, one may first obtain cumulants  and then moments. 

Since for noncommutative random variables, we start with their moments as we have seen so far, it is very important to investigate moment and cumulant sequences for further calculations. Before going to more details, let us see some basic concepts about partitions of an index 
set, which plays an important role in free probability theory.

\subsection{Partitions, Non-crossing Partitions, and Free-Cumulants}

For a positive integer $n$, we denote $[n]\overset{\Delta}{=}\{1,2,...,n\}$. A partition $\pi$ of 
set $[n]$ means $\pi=\{V_1,...,V_k\}$ such that $V_1,...,V_k\subset [n]$ with
$V_i\neq \emptyset$, $V_i\cap V_j=\emptyset$ for all $1\leq i\neq j\leq n$, and 
$V_1\cup \cdots \cup V_k=[n]$. Subsets $V_1,...,V_k$ are called the blocks of $\pi$ and 
$\#(\pi)$ denotes the number of the blocks of $\pi$. ${\cal P}(n)$
denotes the set of all the partitions of $[n]$. 
A partition is called 
a pairing if its each block has size $2$ and the set of all the pairings of $[n]$ is denoted by ${\cal P}_2(n)$.

Let $\pi\in {\cal P}(n)$ and $\{k_i\}_i$ be a sequence. We denote 
$k_{\pi}=k_1^{r_1}k_2^{r_2}\cdots k_n^{r_n}$ where $r_i$ is the number of blocks 
of $\pi$ of size $i$. Then, the determination formulas in (\ref{3.4})-(\ref{3.5}) of moments and cumulants can be re-formulated as
\begin{eqnarray}
 m_n & = & \sum_{\pi\in {\cal P}(n)} k_{\pi}, \label{3.6}\\
 k_n & = & \sum_{\pi \in {\cal P}(n)} (-1)^{\#(\pi)-1}(\#(\pi)-1)! m_{\pi}. \label{3.7}
\end{eqnarray}

For $\pi\in {\cal P}(n)$, denote the moment of $n$ random variables
$x_1,...,x_n$ with partition $\pi$ as 
\begin{equation}\label{3.8}
E_{\pi}(x_1,...,x_n)\overset{\Delta}{=}\mathop{\prod_{V\in {\pi}}}_{V=(i_1,...,i_l)} E(x_{i_1}\cdots x_{i_l}),
\end{equation}
where $V=(i_1,...,i_l)$ means that set $V$ has $l$ distinct elements 
with increasing order as $i_1<i_2<\cdots<i_l$.

When $\pi\in {\cal P}_2(2k)$, i.e., $\pi$ is a pairing of $[2k]$, we have 
\begin{equation}\label{3.9}
E_{\pi}(x_1,..., x_{2k})=\prod_{(i,j)\in \pi}E(x_ix_j).
\end{equation}
With this notation, for Gaussian random variables $X_1,X_2,...,X_n$, we have the following Wick's formula:
\begin{equation}\label{3.10}
E(X_{i_1}\cdots X_{i_{2k}})=\sum_{\pi\in {\cal P}_2(2k)} E_{\pi}(X_{i_1},...,X_{i_{2k}}),
\end{equation}
where $i_1,...,i_{2k}\in [n]$. 

Let $\pi\in {\cal P}(n)$. If there exist $i<j<k<l$ such that 
$i$ and $k$ are in one block $V$ of $\pi$, and $j$ and $l$ in another 
block $W$ of $\pi$, we call that $V$ and $W$ cross. If one cannot find 
any pair of blocks in $\pi$ that cross, partition $\pi$ is called non-crossing.
Denote the set of all non-crossing partitions of $[n]$ by $NC(n)$ and the set
of all non-crossing pairings of $[n]$ by $NC_2(n)$. 

The partition set ${\cal P}(n)$ of $[n]$ is partially ordered via
$$
\pi_1\leq \pi_2 \mbox{ if and only if each block of }\pi_1 \mbox{ is contained in a block of }\pi_2.
$$
With this order, $NC(n)$, as a subset of ${\cal P}(n)$, is 
also partially ordered. The largest and the smallest partitions in both  ${\cal P}(n)$ and $NC(n)$ 
are $[n]$ and $\{ \{1\},\{2\},...,\{n\}\}$,  denoted as $1_n$ and $0_n$, respectively.

\begin{DD}\label{def3}
The following free cumulants $\kappa_n(x_1,...,x_n)$ are defined inductively in terms of moments by the moment-cumulant formula:
\begin{equation}\label{3.11}
E(x_1\cdots x_n)=\sum_{\pi\in NC(n)} \kappa_{\pi}(x_1,...,x_n),
\end{equation}
where 
\begin{equation}\label{3.12}
\kappa_{\pi}(x_1,...,x_n) \overset{\Delta}{=} \mathop{\prod_{V\in \pi}}_{V=(i_1,...,i_l)}
  \kappa_l(x_{i_1},...,x_{i_l}).
\end{equation}
\end{DD}

The above inductive definition is not hard to implement as follows.

For $n=1$, we have $E(x_1)=\kappa_1(x_1)$. Thus, $\kappa_1(x_1)=E(x_1)$.

For $n=2$, we have 
$$
E(x_1x_2)=\kappa_{(1,2)}(x_1,x_2)+\kappa_{(1),(2)}(x_1,x_2)=\kappa_2(x_1,x_2)
+\kappa_1(x_1)\kappa_1(x_2).
$$
Thus,
$$
\kappa_2(x_1,x_2)=E(x_1x_2)-E(x_1)E(x_2), 
$$
etc.

Let $\mu(\pi_1,\pi_2)$ be the M\"{o}bius function on ${\cal P}(n)$ 
\cite{spe4,spe5,hia} that has a recursion formula to calculate. 
Then, we also have the following M\"{o}bius inversion formula: 
\begin{equation}\label{3.13}
\kappa_n(x_1,...,x_n)=
\sum_{\pi\in NC(n)}\mu(\pi, 1_n) E_{\pi}(x_1,...,x_n).
\end{equation}

The moment-cumulant formulas (\ref{3.11}) and (\ref{3.13}) 
for momemts and free-cumulants for noncummtative random variables 
 are 
in analogous to (\ref{3.6}) and (\ref{3.7}) (or
(\ref{3.4}) and (\ref{3.5})) for classical random variables 
in classical probability theory. 

\begin{TT}\label{thm2}
Random variables $x_1,...,x_n$ are free if and only if all mixed 
cumulants of $x_1,...,x_n$  vanish. In other words, $x_1,...,x_n$ are free
if and only if,  for any $i_1,...,i_p\in [n]=\{1,2,...,n\}$ 
with $i_j\neq i_l$ for some $j,l\in [p]$, we have 
$\kappa_p(x_{i_1},...,x_{i_p})=0$.
\end{TT}

The result in the above theorem significantly simplifies the calculations of
the free cumulants of multiple free random variables and therefore, helps to 
calculate the joint moments of multiple free random variables.
For example, if $x$ and $y$ are free, then we have
\begin{eqnarray}
\kappa_n^{x+y} & \overset{\Delta}{=} &
\kappa_n(x+y,...,x+y) \nonumber\\
 & = & \kappa_n(x,...,x)+\kappa_n(y,...,y)+(\mbox{mixed cumulants in }x, y) \nonumber\\
 & = & \kappa_n^x+ \kappa_n^y. \label{3.133}
\end{eqnarray}

\begin{DD}\label{def4}
Let $I$ be an index set. A self-adjoint family $(s_i)_{i\in I}$ is called 
a semicircular family of covariance matrix $C=(c_{ij})_{i,j\in I}$ if $C$ is
non-negative definite and for any $n\geq 1$ and any $n$-tuple $i_1,...,i_n\in I$ we have 
\begin{equation}\label{3.14}
E(s_{i_1}\cdots s_{i_n})=\sum_{\pi\in NC_2(n)} E_{\pi}(s_{i_1},...,s_{i_n}),
\end{equation}
where 
\begin{equation}\label{3.15}
E_{\pi}(s_{i_1},...,s_{i_n})=\prod_{(p,q)\in \pi} c_{i_p,i_q}.
\end{equation}
If $C$ is diagonal, then $(s_i)_{i\in I}$ is a free semicircular family.
\end{DD}

The above formula is the free analogue of Wick's formula for Gaussian 
random variables. If we let $X_1,...,X_r$ be $N\times N$ matrices 
of all entries in all matrices i.i.d. Gaussian random variables, then
they jointly converge in distribution to a free semi-circular family $s_1,...,s_r$ 
of covariance matrix  $(c_{ij})_{1\leq i,j\leq r}=\mbox{I}_r$ where $\mbox{I}_r$
is the identity matrix of size $r$, as $N$ goes to infinity. 
More details on random 
matrices will be seen in Section \ref{sec4}.

\subsection{Cauchy Transforms and R-Transforms}

As we have seen earlier, for classical random variables,
their distributions or density functions can be determined by 
their moment sequences or cumulant sequences as shown in (\ref{3.2}) 
and (\ref{3.3}). To further study noncommutative random variables, 
their moment and cumulant sequences similarly lead to their analytic forms 
as follows.

Let $x$ be a noncommutative random variable and $m_n^x=E(x^n)$ and
 $\kappa_n^x$ be its moments and free cumulants, respectively. 
Their power series (moment and cumulant generating functions)
in an indeterminate $z$ are defined by
\begin{equation}\label{3.16}
M(z)=1+\sum_{n=1}^{\infty} m_n^x z^n  \mbox{  and  }
C(z)=1+\sum_{n=1}^{\infty} \kappa_n^x z^n.
\end{equation}
Then, the following identity holds:
\begin{equation}\label{3.17}
M(z)=C(zM(z)).
\end{equation}

The Cauchy transform of $x$ is defined by
\begin{equation}\label{3.18}
G(z)\overset{\Delta}{=} E\left( \frac{1}{z-x}\right)
=\sum_{n=0}^{\infty}\frac{E(x^n)}{z^{n+1}}=\sum_{n=0}^{\infty}\frac{m_n^x}{z^{n+1}}=z^{-1}M(z^{-1}),
\end{equation}
and the $R$-transform of $x$ is defined by
\begin{equation}\label{3.19}
R(z)\overset{\Delta}{=}\frac{C(z)-1}{z}=\sum_{n=0}^{\infty} \kappa_{n+1}^x z^n.
\end{equation}
If we let $K(z)\overset{\Delta}{=}R(z)+z^{-1}$, then $K(G(z))=z$, i.e., $K(z)$ is the inverse of the Cauchy transform $G(z)$. 

If we let $G_x(z)$ and $R_x(z)$ denote the Cauchy transform
and the $R$-transform  of random variable $x$, respectively, then,
for two free random variables $x$ and $y$, from (\ref{3.133}) we have 
\begin{equation}\label{3.20}
R_{x+y}(z)= R_x(z)+R_y(z).
\end{equation}

In case not both $R_x(z)$ and $R_y(z)$ are well-defined on a region of $z$, 
one may be able to 
find the Cauchy transform $G_{x+y}(z)$ of $x+y$ for free random variables $x$ and $y$ 
from the Cauchy transforms $G_x(z)$ and $G_y(z)$ of $x$ and $y$ as follows.

We shall see soon below that when $z$ is in the upper complex plane 
$\mathbb{C}^+ \overset{\Delta}{=} \{ c\in \mathbb{C}| \mbox{Im}(c)>0\}$ 
where $\mathbb{C}$ stands
for the complex plane and Im stands for the imaginary part of a complex numnber, a Cauchy transform is well-defined.

For an $z\in \mathbb{C}^+$, solve the following system of two equations 
for two unknown functions $\omega_x(z)$ and $\omega_y(z)$:
\begin{equation}\label{3.21}
G_x(\omega_x(z))=G_y(\omega_y(z)) \mbox{ and } 
\omega_x(z) +\omega_y(z)-\frac{1}{G_x(\omega_x(z))}  =  z.
\end{equation}
Then, 
\begin{equation}\label{3.22}
G_{x+y}(z)=G_x(\omega_x(z))=G_y(\omega_y(z)).
\end{equation}

If noncommutative random variable $x$ is self-adjoint, 
then it has a spectral measure 
$\nu$ on $\mathbb{R}$ such that the moments of $x$ are
the same as the conventional moments of the probability measure 
$\nu$. One can simply see it when $x$ is a Hermitian matrix and then $x$ can be diagonalized by a unitary matrix and has real-valued eigenvalues. These real-valued eigenvalues are the spectra of $x$ that are discrete for a finite matrix but may become continuous when $x$ is a general operator over an 
infinite dimensional space. In this case, we say that random variable $x$ has 
distribution $\nu$. 

Then, the Cauchy transform $G(z)$ of $x$ can be formulated as
\begin{equation}\label{3.23}
G(z)=\int_{\mathbb{R}} \frac{1}{z-t} d\nu(t),
\end{equation}
and $G(z)$ is also called the Cauchy transform of $\nu$. 

One can clearly see 
from (\ref{3.23}) that Cauchy transform $G(z)$ is well-defined when
$z\in \mathbb{C}^+$. In fact, $G(z)$ is analytic in $\mathbb{C}^+$, i.e.,
it exists  derivatives of all orders for any $z\in \mathbb{C}^+$. 
Furthermore, $G(z)\in \mathbb{C}^{-}$, the lower complex plane similarly defined as $\mathbb{C}^+$. 
In other words, a Cauchy transform $G(z)$ maps $\mathbb{C}^+$ to $\mathbb{C}^{-}$.

From (\ref{3.23}), one can also see that the Cauchy transform excludes the real axis
$\mathbb{R}$ for $z$, which is because when $z\in \mathbb{R}$, the integration 
may not exist. After saying so, it may exist in the generalized function sense
as if $z\in \mathbb{R}$, the Cauchy transform (\ref{3.23}) becomes 
the Hilbert transform of $d\nu(t)/dt$. 

When probability measure $\nu$ is compactly supported, i.e., it is supported on a finite interval, not only its Cauchy transform is analytic in $\mathbb{C}^+$, but also its $R$-transform is  analytic on some disk centered at the origin.
This, however, may not be true for a general probability measure $\nu$.  
For more details, see \cite{spe3}.

With a Cauchy transform $G(z)$, its corresponding probability measure can be 
formulated by the Stieltjes inversion formula as follows.

\begin{TT}\label{thm3}
Let $\nu$ be a probability measure on $\mathbb{R}$ and $G(z)$ be its
Cauchy transform. For $a<b$, we have
\begin{equation}\label{3.24}
- \lim_{\tau\rightarrow 0^+}\frac{1}{\pi} \int_a^b \mbox{Im}(G(t+j\tau))dt 
=\nu( (a,b))+\frac{1}{2}\nu(\{a,b\}),
\end{equation}
where $\nu((a,b))$ and $\nu(\{a,b\})$ are the continuous and the discrete 
parts of the measure $\nu$, respectively. 
If $\nu_1$ and $\nu_2$ are two probability measures on $\mathbb{R}$ with
equal Cauchy transforms, i.e., $G_{\nu_1}(z)= G_{\nu_2}(z)$, then $\nu_1=\nu_2$.
\end{TT}

This result tells us that Cauchy transforms and probability measures 
 (distributions or random variables) are 
one-to-one corresponding to each other. 

If $x$ and $y$ are two free self-adjoint random variables with distributions 
$\nu_x$ and $\nu_y$, respectively. The distribution of $x+y$ is called 
the free convolution of those of $x$ and $y$, which is 
denoted by $\nu_x\boxplus \nu_y$. 

As an example of Cauchy transform, when $\nu$ is semicircular with density 
function $q(t)$ in (\ref{2.6}), its Cauchy transform is, \cite{spe3},
\begin{equation}\label{3.25}
G_s(z)=\frac{z-\sqrt{z^2-4}}{2}.
\end{equation}

\section{Application in Random Matrices}\label{sec4}
As mentioned in Introduction, random matrices with entries of complex 
Gaussian random variables are often used in wireless communications and signal 
processing. In particular, their singular value (eigenvalue) distributions 
play an important role in analyzing  wireless communications systems. 
This section is on applying free probability theory to random matrices of 
large sizes. It tells us how to use the second order
statistics of the entries of  random matries to calculate 
their asymptotic eigenvalue distributions. 

\subsection{GUE Random Matrices and Wigner's Semi-Circle Law}

Let $X_N$ be an $N\times N$ matrix with complex random variables 
$a_{ij}=x_{ij}+{\bf i}y_{ij}$ as
entries such that $x_{ij}$ and $y_{ij}$ are real Gaussian random variables,
 $\sqrt{N} a_{ij}$ is a standard complex random variable, i.e.,
$E(a_{ij})=0$ and $E(|a_{ij}|^2)=1/N$ and 
\begin{itemize}
\item[1)] $a_{ij}=a_{ji}^*$,
\item[2)] $\{x_{ij}\}_{i\geq j}\cup \{y_{ij}\}_{i>j}$ are i.i.d.
\end{itemize}
In this case,  $X_N$ is Hermitian, i.e., self-adjoint. $X_N$ is called a Gaussian 
unitary ensemble (GUE) random matrix.
The following theorem is Wigner's semi-circle law.

\begin{TT}\label{thm4}
  If $\{X_N\}_N$ is a sequence of GUE random matrices, then, for any positive 
integer $k$, 
\begin{eqnarray*}
\lim_{N\rightarrow \infty} E(tr(X_N^k)) & = & \frac{1}{2\pi} 
\int_{-2}^2 t^k \sqrt{4-t^2}dt \\
 & = & \left\{ \begin{array}{ll} \frac{1}{l+1} 
{2l \choose l}, & \mbox{if }k=2l \mbox{ for some positive integer }l,\\
 0, & \mbox{if }k \mbox{ is odd}.
\end{array} 
\right.
\end{eqnarray*}
where $tr$ stands for the normalized matrix trace, i.e., 
$tr(\cdot)\overset{\Delta}{=}\mbox{Tr}(\cdot)/N$ with the conventional 
matrix trace Tr. 
\end{TT}

Since $X_N$ is Hermitian, it has spectra (eigenvalues) $\nu_N$
that is a random variable as well.  Since  $tr(X_N^k)=tr(\nu_N^k)$,
we have 
$$
\lim_{N\rightarrow \infty} E(tr(X_N^k))= \lim_{N\rightarrow \infty}
\int_{\mathbb{R}} t^kd\nu_N(t).
$$
Thus, the above theorem says that 
the eigenvalues of $X_N$ converge in distribution
 to the semicircular random variable.
 In fact, the convergence in distribution 
can be made stronger to the almost surely convergence. 

\subsection{Asymptotic Freeness of GUE Random Matrices}

For random matrices $X$ as noncommuntative random variables, their linear functional $E$ used in Section \ref{sec2} is defined as $E(tr(X))$, i.e., $E(\cdot)$ used before for a noncommutative  random variable $x$ 
corresponds to $E(tr(\cdot))$ for a random matrix $X$ in what follows.

\begin{DD}\label{def5}
Let $(X_N)_N$ and $(Y_N)_N$ be two sequences of $N\times N$ matrices.
We say that $X_N$ and $Y_N$ are  asymptotically free 
if they converge in distribution to two free random variables $x$ and $y$,
respectively, as $N$ goes to infinity. 
\end{DD}

From Definitions \ref{def2} and \ref{def5}, $X_N$ and $Y_N$ 
are asymptotically free, if for any positive integer $m$
and non-negative integers $p_1,q_1,...,p_m,q_m$ we have
$$
\lim_{N\rightarrow \infty} E(tr(X_N^{p_1}Y_N^{q_1}\cdots X_N^{p_m}Y_N^{q_m}))
=E(x^{p_1}y^{q_1}\cdots x^{p_m}y^{q_m}),
$$
for two free random variables $x$ and $y$. 

For a sequence of $N\times N$ deterministic matrices $(D_N)_N$,
if  $\lim_{N\rightarrow \infty}tr(D_N^m)$ exists for every non-negative 
integer $m$, we say $D_N$ converges to $d$ 
in distribution, where $d$ is a noncommutative random variable 
and its $m$th moment is the same as the limit. We also write it as 
$\lim_{N\rightarrow \infty} D_N \overset{distr}{=} d$
or $D_N \overset{distr,}{\longrightarrow} d$. 

With the above notations, the following theorem of Voiculescu improves Wigner's semi-circle law.

\begin{TT}\label{thm5}
Assume $X_N^{(1)},...,X_N^{(p)}$ are $p$ independent $N\times N$ GUE 
random matrices and $D_N^{(1)},...,D_N^{(q)}$ are $q$ 
deterministic $N\times N$ matrices such that
$$
D_N^{(1)},...,D_N^{(q)} \overset{distr}{\longrightarrow} d_1,...,d_q \mbox{ as }
N\rightarrow \infty.
$$
Then,
$$
X_N^{(1)},...,X_N^{(p)}, D_N^{(1)},...,D_N^{(q)} \overset{distr}{\longrightarrow} 
s_1,...,s_p, d_1,...,d_q \mbox{ as }
N\rightarrow \infty,
$$
where each $s_i$ is semicircular and $s_1,...,s_p, \{d_1,...,d_q\}$ are free.
The convergence above also holds almost surely. 
\end{TT}

This result tells that independent GUE random matrices $X_N^{(1)},...,X_N^{(p)},
\{D_N^{(1)},...,D_N^{(q)}\}$ are asymptotically free
when $N$ is large. Furthermore, $X_N^{(1)},...,X_N^{(p)}$
asymptotically have the same distributions as free semicircular 
elements $s_1,...,s_p$ do, and 
this is still true even when they are mixed with deterministic matrices.

\subsection{Asymptotic Freeness of Haar Distributed Unitary Random Matrices}

For a general Hermtian random matrix, it can be diagonalized by a unitary 
matrix and in this case, the unitary matrix is  random as well. 
Therefore, it is  also important to study unitary random matrices. 

Let ${\cal U}(N)$ denote the group of $N\times N$ unitary matrices $U$, i.e.,
$UU^*=U^*U=\mbox{I}_N$. Since ${\cal U}(N)$ is bounded (compact), it has Haar meansure 
 $dU$ with $\int_{{\cal U}(N)} dU=1$. Thus,  $dU$ is a probability measure (it can be understood as a uniform distribution). A Haar distributed 
unitary random matrix is a matrix $U_N$ randomly chosen in ${\cal U}(N)$
with respect to Haar measure. One method to construct Haar unitary matrices is
as follows. First, take an $N\times N$ random matrix whose entries are the 
independent standard complex Gaussian random variables. Then, use the Gram-Schmidt orthogonalization procedure to make it unitary. 

A noncommutative random variable $u$ is called Haar unitary if 
it is unitary, i.e., $uu^*=u^*u=1$ and $E(u^m)=\delta_{0,m}$, i.e., $0$ when $m>0$. A Haar unitary random matrix is Haar unitary, i.e., if $U\in {\cal U}(N)$,
then $E(tr(U^m))=0$ for $m>0$ \cite{spe3}.

\begin{TT}\label{thm6}
Assume $U_N^{(1)},...,U_N^{(p)}$ are $p$ independent $N\times N$ Haar unitary  
random matrices and $D_N^{(1)},...,D_N^{(q)}$ are $q$ 
deterministic $N\times N$ matrices such that
$$
D_N^{(1)},...,D_N^{(q)} \overset{distr}{\longrightarrow} d_1,...,d_q \mbox{ as }
N\rightarrow \infty.
$$
Then, as $N \rightarrow \infty$, 
$$
U_N^{(1)},U_N^{(1)*},...,U_N^{(p)},U_N^{(p)*},  D_N^{(1)},...,D_N^{(q)} \overset{distr}{\longrightarrow} 
u_1,u_1^*,...,u_p,u_p^*, d_1,...,d_q, 
$$
where each $u_i$ is  Haar unitary and $\{u_1,u_1^*\},...,\{u_p,u_p^*\}, \{d_1,...,d_q\}$ are free.
The convergence above also holds almost surely. 
\end{TT}

A more special case is as follows.

\begin{TT}\label{thm7}
Let $A_N$ and $B_N$ be two sequences of deterministic $N\times N$ 
matrices with $\lim_{N\rightarrow \infty}A_N \overset{distr}{=}a$
and $\lim_{N\rightarrow \infty}B_N \overset{distr}{=}b$.
Let $U_N$ be a sequence of $N\times N$ Haar unitary random matrices.
Then, 
$$
A_N, U_NB_NU_N^* \overset{distr}{\longrightarrow} a,b \mbox{ as } N \rightarrow \infty, 
$$
where $a$ and $b$ are free. 
This convergence also holds almost surely.
\end{TT}

The above theorem says that $A_N$ and $U_NB_NU_N^*$ are
asymptotically free when $N$ is large. 

\subsection{Aymptotic Freeness of Wigner Random Matrices}

Let $\mu$ be a probability measure on $\mathbb{R}$ and $a_{ij}$ with $i\leq j$ 
be i.i.d. real random variables with distribution $\mu$. Let $a_{ij}=a_{ji}$
for $i>j$, and 
$$A_N=\frac{1}{\sqrt{N}} (a_{ij})_{1\leq i,j\leq N},
$$
which is self-adjoint (symmetry) and  called Wigner random matrix (ensemble).

\begin{TT}\label{thm8}
Let $\mu_1,...,\mu_p$ be probability 
measures on $\mathbb{R}$ with all moments exist and $0$ mean. 
Assume $A_N^{(1)},...,A_N^{(p)}$ are $p$ independent $N\times N$ Wigner 
random matrices with entry distributions $\mu_1,...,\mu_p$,
respectively,  and $D_N^{(1)},...,D_N^{(q)}$ are $q$ 
deterministic $N\times N$ matrices such that
$$
D_N^{(1)},...,D_N^{(q)} \overset{distr}{\longrightarrow} d_1,...,d_q \mbox{ as }
N\rightarrow \infty,
$$
and 
$$
\sup_{r,N}\|D_N^{(r)}\|<\infty.
$$
Then, as $N \rightarrow \infty$, 
$$
A_N^{(1)},...,A_N^{(p)}, D_N^{(1)},...,D_N^{(q)} \overset{distr}{\longrightarrow} 
s_1,...,s_p,, d_1,...,d_q, 
$$
where each $s_i$ is semicircular and $s_1,...,s_p, \{d_1,...,d_q\}$ are free.
\end{TT}

As a special case, $A_N D_N A_N, E_N \overset{distr}{\longrightarrow} sds, e$,
where $s$ is semicircular, $sds$ and $e$ are free, and $e$ can be arbitrary.

\section{Free Deterministic Equivalents and Random Matrix Singular Value Distribution Calculations}\label{sec5}

Let $H$ be an $N\times M$ wireless channel matrix, which is usually modelled
as a random matrix, with additive white Gaussian noise (AWGN) of  variance 
$\sigma$. Then, its mutual information is
\begin{equation}\label{5.1}
C(\sigma)=\frac{1}{N}E\left[\log \det \left( \mbox{I}_N+ \frac{HH^*}{\sigma}\right) \right],
\end{equation}
where $^*$ stands for Hermitian operation. Let $\nu(\lambda)$ denote
the eigenvalue distribution (or spectra, or probability measure) of
matrix $HH^*$. Then, when $N$ is large, 
\begin{equation}\label{5.2}
C(\sigma)=\int_{0}^{\infty} \log\left( 1+\frac{\lambda}{\sigma}\right)
d\nu(\lambda).
\end{equation}
On the other hand, the Cauchy transform of the probability measure 
$\nu$ and matrix $HH^*$ is
\begin{equation}\label{5.3}
G(z)=\int_{0}^{\infty} \frac{1}{z-\lambda} d\nu(\lambda)
=E(tr(z\mbox{I}_N-HH^*)^{-1}),
\end{equation}
where $z\in \mathbb{C}^+$. Assume that $G(z)$ exists as $\mbox{Im}(z)\rightarrow 0^+$, whose limit is denoted by $G(\omega)$ with $\omega=$Re$(z)$. 
For semicircular distribution, from (\ref{3.25}) one can see that
$G(\omega)$ exists when $\omega=$Re$(z)>2$. 
Then, \cite{tul2}, 
\begin{equation}\label{5.4}
C(\sigma)=\int_{\sigma}^{\infty}\left( \frac{1}{\omega}-G(-\omega)\right)
d\omega.
\end{equation}
The above  formula tells us that, to calculate the mutual information 
of the channel with channel matrix $H$, we only need to calculate the 
Cauchy transform of matrix $HH^*$. 

As an example, if $HH^*$ is a GUE random matrix, then, when $N$ is large,
it is approximately  semicircular and its Cauchy transform has the form of 
(\ref{3.25}) with a proper normalization. Thus, its mutual information 
can be calculated. However, in applications, $HH^*$ may not be 
 a GUE matrix. We next introduce free deterministic equivalents
to help to calculate  the Cauchy transforms of large random matrices,
such as the above $HH^*$, 
 based
on Speicher \cite{spe1, spe3, spe6}. 

\subsection{Matrix-Wise Free Deterministic Equivalents}
From Section \ref{sec4}, we know that when the entries $X_{ij}$, $i\geq j$, 
of an $N\times N$
self-adjoint (symmetry for real-valued or Hermitian for complex-valued)
matrix $X$ are i.i.d. random variables, when $N$ is large, it is approximately
semicircular. It is also true for multiple such random 
matrices  and multiple deterministic matrices jointly. 

For a non-adjoint random matrix $X$ of i.i.d. Gaussian entries, 
it can be made into two independent 
self-adjoint GUE matrices as $Y_1=(X+X^*)/\sqrt{2}$
and $Y_2=-{\bf i}(X-X^*)/\sqrt{2}$. 
Then, $X=(Y_1+{\bf i}Y_2)/\sqrt{2}$. 
In this case, $X$ converges in distribution to 
$s=(s_1+{\bf i} s_2)/\sqrt{2}$ for two free semicircular elements 
$s_1$ and $s_2$  with the same distribution. While $s_1$ and $s_2$ are 
semicircular, we call $s$ circular.

In \cite{spe6,spe3} it is proposed to replace these random matrices 
by semicircular and circular elements etc. 
Consider the following collections of $N\times N$ matrices,
where for each random matrix, its entries of different random 
variables are i.i.d.:
\begin{eqnarray*}
& {\bf X}=\{X_1,...,X_{n_1}\}  : & \mbox{independent self-adjoint  matrices},\\
& {\bf Y}=\{Y_1,...,Y_{n_2}\}  : & \mbox{independent non-self-adjoint  matrices},\\
& {\bf U}=\{U_1,...,U_{n_3}\}  : & \mbox{independent Haar distribued unitary  matrices},\\
& {\bf D}=\{D_1,...,D_{n_4}\}  : & \mbox{deterministic  matrices}.
\end{eqnarray*}
Let 
\begin{eqnarray*}
{\bf s}=\{s_1,...,s_{n_1}\} &: & \mbox{free semicircular},\\
{\bf c}=\{c_1,...,c_{n_2}\} &: & \mbox{free circular},\\
{\bf u}=\{u_1,...,u_{n_3}\} &: & \mbox{free Haar unitary},\\
{\bf d}=\{d_1,...,d_{n_4}\} &: & \mbox{abstract elements},
\end{eqnarray*}
Assume that the joint distribution of ${\bf D}$ is the same as that of 
${\bf d}$, and ${\bf X},{\bf Y},{\bf U}$ are independent among each other.
Also assume that ${\bf s},{\bf c}, {\bf u}$ have their each individual 
distribution 
asymptotically 
the same as that of ${\bf X}, {\bf Y}, {\bf U}$, respectively.

Let $P_N$ be a multi-variable polynomial of ${\bf X}, {\bf Y}, {\bf U}, {\bf D}$. Then, when $N$ is large, 
$$P_N=P(X_1,...,X_{n_1},Y_1,...,Y_{n_2},U_1,...,U_{n_3},D_1,...,D_{n_4})$$
can be replaced by
$$
P_N^{\Box}=P(s_1,...,s_{n_1},c_1,...,c_{n_2},u_1,...,u_{n_3},d_1,...,d_{n_4})$$
and 
$P_N^{\Box}$ is called the (matrix-wise) free deterministic equivalent of $P_N$. Then, we have, for any positive integer $k$, 
$$
\lim_{N\rightarrow \infty} E(tr(P_N^k))=E((P_N^{\Box})^k).
$$

Now let us go back to the matrix $HH^*$ in (\ref{5.1}). Although 
matrix $H$ is not self-adjoint itself, but if we follow
\cite{spe1} and \cite{spe3} and let 
\begin{equation}\label{5.5}
T=\left( \begin{array}{cc} 0 & H \\ H^* & 0 \end{array} \right),
\end{equation}
then, matrix $T$  is self-adjoint. Furthermore, 
\begin{equation}\label{5.6}
T^2=\left( \begin{array}{cc} HH^* & 0 \\ 0 & H^*H \end{array} \right),
\end{equation}
which includes $HH^*$ as a diagonal block.
 Using  operator-valued 
free probability theory \cite{spe6,spe3}, it can be similarly treated 
as what is done in the previous sections. 
Note that $T^2$ is just a polynomial of $T$ but 
unfortunately not all entries in matrix $T$ 
have the same distribution, which makes the above matrix-wise 
free deterministic equivalent approach difficult to use.
In order to deal with this problem, we next consider component-wise
free deterministic equivalents.

\subsection{Component-Wise Free Deterministic Equivalents and 
Cauchy Transform Calculation of Random Matrices}

This part is mainly from \cite{spe1}. 
We consider  $N\times N$ random matrices $X=(X_{ij})$
where $X_{ij}$ are  complex Gaussian random variables
with $E(X_{ij})=0$ and $E(X_{ij}X_{ij}^*)=\sigma_{ij}/N$,
where $\sigma_{ij}$ are independent of $N$. 
 Now we replace all entries $X_{ij}$ in $X$ by
(semi)circular elements $c_{ij}$ such that 
\begin{equation}\label{5.7}
E(c_{ij}c_{ij}^*)=E(X_{ij}X_{ij}^*)=\sigma_{ij}/N
\end{equation}
where if $X_{ij}$ is real-valued (or complex-valued), then, 
$c_{ij}$ is semicircular (or circular) with mean $0$; 
if $X_{ij}$ and $X_{kl}$ are independent, then $c_{ij}$ and $c_{kl}$
are free; if $X_{ij}=X_{kl}$, then $c_{ij}=c_{kl}$; and 
$E(X_{ij}X_{kl}^*)=E(c_{ij}c_{kl}^*)=(E(c_{ij}^*c_{kl}))^*$. 
Then,  we form 
an $N\times N$
 matrix of (semi)circular elements as $c=(c_{ij})$.  Matrix 
$c$ is called the component-wise free deterministic equivalent of 
matrix $X$. 

Let $X_1,...,X_{n_1}$ be $n_1$ random matrices, each of which is 
 specified above,  where all entries of each of these matrices  are 
independent from all entries of all the remaining matrices.
Let $c_1,...,c_{n_1}$ be the component-wise deterministic equivalents 
of $X_1,...,X_{n_1}$, and all elements in $c_i$ are free from all elements 
in $c_j$ when $i\neq j$. 
Let $D_1,...,D_{n_2}$ be $n_2$ deterministic matrices.
Assume $P_N$ is a multi-variable polynomial and
\begin{eqnarray}
P_N & = & P(X_1,...,X_{n_1},D_1,...,D_{n_2}),  \label{5.8}\\
P_N^{\Box} & = & P(c_1,...,c_{n_1}, D_1,...,D_{n_2}). \label{5.9}
\end{eqnarray}
We call that $P_N^{\Box}$ is the component-wise free deterministic 
equivalent of $P_N$. 

\subsubsection{Independent Cases}

Consider the case when every matrix $X_i$ is Hermitian/self-adjoint
 and entries $X_{ij}$ for $i\geq j$ are all independent. 
It is explicitly  shown in \cite{anan} that 
$\lim_{N\rightarrow \infty} (P_N-P_N^{\Box})=0$, i.e.,
the  matrices 
$X_1,...,X_{n_1}, D_1,...,D_{n_2}$ have the same joint 
distribution as matrices $c_1,...,c_{n_1},D_1,...,D_{n_2}$ do. 
Thus, $c_1,...,c_{n_1}$ may be used to calculate 
the Cauchy transforms of $X_1,...,X_{n_1}$ when only 
the variances of the entries in matrices $X_i$ are used, 
as $N$ is large. 

We now consider a special example shown in \cite{spe1}. 
Let $X=(X_{ij})$ be an $N\times N$ Hermitian/self-adjoint
 Gaussian random matrix 
with $E(X_{ij})=0$ and $E(X_{ij}X_{ij}^*)=\sigma_{ij}/N$, and 
let $c=(c_{ij})$ be its component-wise deterministic 
equivalent, i.e., $E(c_{ij})=0$ and $E(c_{ij}c_{ij}^*)=\sigma_{ij}/N$.
Note that since $X_{ij}=X_{ji}^*$, we have $c_{ij}=c_{ji}^*$ as well.
Let $A$ be an $N\times N$ deterministic matrix.

Consider the matrix sum  $Y=A+X$. We next show how to calculate the Cauchy 
transform of $Y$ by calculating that of $T=A+c$. 

For an $N\times N$  deterministic matrix $B=(B_{ij})$, define a mapping
$\eta$ that maps $B$ to another $N\times N$ deterministic matrix $\eta(B)$
with its $(i,j)$th component as
\begin{equation}\label{5.10}
[\eta(B)]_{ij}\overset{\Delta}{=}E(cBc)=\sum_{k,l}E(c_{ik}B_{kl}c_{lj})
=\sum_{k,l} E(c_{ik}c_{jl}^*)B_{kl}
=\delta_{i,j}\sum_k \sigma_{ik}B_{kk},
\end{equation}
which shows that $\eta(B)$ is a diagonal matrix. 
Then, the Cauchy transform $g_T(z)$ of $T$ can be determined
by solving the following fixed point equation \cite{spe7}, \cite{spe1}: 
\begin{eqnarray}
g_T(z) & = & tr(G_T(z)), \label{5.11}\\
G_T(z) & = & E\left( \frac{1}{z-\eta(G_T(z))-A}\right), \label{5.12}
\end{eqnarray}
where $E(B)\overset{\Delta}{=}(E(B_{ij}))$. It is shown in
\cite{spe7} that there is exactly one solution of the above 
fixed point equation with the proper positivity constriant.

We next consider the case when $X$ is not Hermitian, such as the channel matrix $H$ in (\ref{5.1}), where all entries of $X$ are independent. In this case, 
consider $Y=A+X$ and we next calculate the Cauchy transform
of $YY^*$. To do so, define 
\begin{equation}\label{5.122}
T=\left( \begin{array}{cc} 0 & Y\\ Y^* & 0 \end{array} \right)
=\left( \begin{array}{cc} 0 & A\\ A^* & 0 \end{array} \right)
+\left( \begin{array}{cc} 0 & X\\ X^* & 0 \end{array} \right).
\end{equation}
Then,
\begin{equation}\label{5.13}
T^2=\left( \begin{array}{cc} YY^* & 0\\ 0 & Y^*Y \end{array} \right).
\end{equation}
Since the eigenvalue distributions of $YY^*$ and $Y^*Y$
are the same, the Cauchy transform of $YY^*$ is the same as that 
of $T^2$. It is presented in \cite{spe1} as follows.

For an $M\times M$ matrix $B=(B_{ij})$, define
$$
E_{D_M}(B)\overset{\Delta}{=} \mbox{diag}(E(B_{11}),\cdots, E(B_{MM})),
$$
where diag stands for  the $M\times M$ diagonal matrix with 
its arguments as its diagonal elements, and also define 
\begin{equation}\label{5.100}
\eta_1(B)\overset{\Delta}{=} E(cBc^*) \mbox{ and } 
\eta_2(B)\overset{\Delta}{=} E(c^*Bc).
\end{equation}
Note that since all the entry elements in matrix $c$ are free from each other,
$E(cBc^*)=E_{D_N}(cBc^*)$ and $E(c^*Bc)=E_{D_N}(c^*Bc)$ as what is shown 
for $\eta$ in (\ref{5.10}). 

Then, the Cauchy transform  $g_{T^2}(z)$ of $T^2$ or $YY^*$ is 
$g_{T^2}(z)=tr(G_{T^2}(z))$ and $G_{T^2}(z)$ can be obtained by solving the 
following fixed point equations \cite{hac}, \cite{spe1}:
\begin{equation}\label{5.144}
zG_{T^2}(z^2) 
=G_T(z) =
E_{D_{2N}} \left[ \left(
\begin{array}{cc} z-z\eta_1(G_2(z^2)) & -A \\
                        -A^* & z-z\eta_2(G_1(z^2))  \end{array}
\right)^{-1} \right],
\end{equation}
where 
\begin{eqnarray}
 zG_1(z) & = & E_{D_N} \left[ \left( 1-\eta_1(G_2(z))+
 A \frac{1} {z-z\eta_2(G_1(z))}A^*
\right)^{-1} \right], \label{5.14}\\
 zG_2(z) & = & E_{D_N} \left[ \left( 1-\eta_2(G_1(z))+
 A \frac{1} {z-z\eta_1(G_2(z))}A^*
\right)^{-1} \right]. \label{5.15}
\end{eqnarray}

\subsubsection{Correlated Cases and Summary}

When the entries in matrix $X$ are correlated, 
similar treatment as the above can be done \cite{spe1}.
One can still 
get the Cauchy transform of $Y$ when $X$ is Hermitain 
by solving the fixed point equation (\ref{5.11})-(\ref{5.12})
and the Cauchy transform of $YY^*$ when $X$ is not Hermitian 
by solving  the fixed point equation (\ref{5.144})-(\ref{5.15}),
where $\eta(B)=E(cBc)$ may not be diagonal as what is calculated 
in (\ref{5.10}), and $\eta_1(B)$ and $\eta_2(B)$ may not be diagonal either.
An example of correlated entries in $X$ is that each column vector (or row
vector) of $X$ is a linear transform of a vector of independent 
Gaussian random variables. 

A simpler example of correlated cases is when random matrix $X_1=BX$ where
$B$ is a deterministic matrix and $X$ is a random matrix of independent 
entries. 
In this case, $X_1$ can be treated as a product of two matrices of $B$ and $X$
 and thus, was covered previously.

The above Cauchy transform calculation is only based on the 
covariances (the second order statistics) of the entries of random matrix $X$.
 As we mentioned easlier,
in this case one does not need to implement Monte-Carlo simulations to 
do the calculations that may be not convenient in practice when 
$X$ has a large size.

Going back to the mutual information in the beginning of this section, 
we can just let $A=0$ in the above to get the Cauchy transform
of $HH^*=YY^*$. 

As a  remark, the deterministic equivalents defined above are from 
\cite{spe6, spe3, spe1}, which we refer to for any difference  with
those appeared in \cite{gir, hac}.

\section{Conclusions}\label{sec6}

As mentioned in the beginning of this paper, the main goal
here is to introduce free probability theory and its application 
to random matrices as simple as possible. It is for a non-mathematics 
major researcher in, for example, communications and signal processing areas.
This paper is mainly based on \cite{spe1, spe2, spe3, spe4}.

Free probability theory is about noncommutative elements or random variables, 
such as, 
random matrices, in contrast to the conventional 
(real-valued or complex-valued) commutative random variables in 
the classical probability theory. The freeness significantly simplifies 
the calculations of the moments and therefore the distributions, 
and interestingly, random matrices,  
when their size is large,  do have the freeness asymptotically. 
Therefore, free probability theory is naturally applied to 
calculate the asymptotic 
distributions of the eigenvalues/singular-values of random matrices when their 
size is large, such as wireless channel matrices in massive MIMO systems.
It is particularly interesting that the calculation only needs the second 
order statistics of the matrix entries. 

This paper is based on the author's own understanding on free 
probability theory and by no means the material covered in this paper 
is complete. More complete materials on this topic are  
 referred to \cite{spe1, spe2, spe3, spe4, spe5, hia, tul2,  cou1}. 

\begin{center}
{\bf Acknowledgement}
\end{center}

The author would like to thank Dr. Roland Speicher for his free online video 
lectures \cite{spe4} and for his helps to my questions 
and  Dr. Anan Lu for his numerous discussions on free probability theory.

\end{document}